\documentclass[11pt]{amsart}
\usepackage{amssymb}

\newtheorem{theorem}{{\bf Theorem}}
\newtheorem{lemma}{{\bf Lemma}}
\newtheorem{prop}{{\bf Proposition}}
\newtheorem{cor}{{\bf Corollary}}
\begin{document}
\title{Twisted conjugacy in free groups and Makanin's question}
\author[Bardakov]{Valerij Bardakov}
\address{Sobolev Institute of Mathematics, Novosibirsk 630090, Russia}
\email{bardakov@math.nsc.ru}
\author[Bokut]{Leonid Bokut}
\address{Sobolev Institute of Mathematics, Novosibirsk 630090, Russia}
\email{bokut@math.nsc.ru}
\author[Vesnin]{Andrei Vesnin}
\address{Sobolev Institute of Mathematics, Novosibirsk 630090,
Russia; and School of Ma\-the\-matical Sciences, Seoul National University,
Seoul 151-747, Korea} \email{vesnin@math.nsc.ru, vesnin@math.snu.ac.kr}
%
%
\thanks{Authors were supported in part by the Russian Foundation
for Basic Research (grants~01--01--00630 and~02--01--01118).}
\date{\today}

\begin{abstract}
We discuss the following question of G.~Makanin from
``Kou\-rov\-ka notebook'': does there exist an algorithm to
determine is for an arbitrary pair of words $U$ and $V$ of a free
group $F_n$ and an arbitrary automorphism $\varphi \in \mbox{\rm
Aut} (F_n)$ the equation $\varphi (X) U = V X$ solvable in $F_n$?
We give the affirmative answer in the case when an automorphism is
virtually inner, i.e. some its non-zero power is an inner
automorphism of $F_n$.
\end{abstract}

\maketitle




\section{Conjugacy and twisted conjugacy}

Suppose $G$ is a group given by a presentation in generators and
defining relations. Three following  decision problems formulated
by  M.~Dehn \cite{Dehn} in 1912 (see also \cite[Ch.~1, \S~2;
Ch.~2, \S~1 ]{LS}) are fundamental in the group theory.

\underline{Word problem:} {\it Does there exist an algorithm to determine if an
arbitrary group word $W$ given in the generators of $G$ defines the identity
element of $G$?}

\underline{Conjugacy problem:} {\it Does there exist an algorithm
to determine is an arbitrary pair of group words $U$, $V$ in the
generators of $G$ define conjugate elements of $G$?}

\underline{Isomorphism problem:} {\it Does there exist an algorithm to
determine for any two arbitrary finite presentations whether the groups they
present are isomorphic or not?}

All three of these problems have negative answers in general (see,
for example \cite{AD}, \cite[Ch.~6.7]{BK}). These results together
with solutions of Dehn's problems in restricted cases have been of
central importance in the combinatorial group theory. For this
reason combinatorial group theory has always searched for and
studied classes of groups in which these decision problems are
solvable. Also, various generalizations of these problems were
considered (see, for example \cite{KhS}, \cite[Ch.~4, \S~4]{LS}).

In the present paper we discuss the twisted conjugacy problem,
that can be con\-si\-dered as a generalization of Dehn's conjugacy
problem. Suppose $G$ is a group given by a presentation and $H$ is
a subset of its automorphisms group $\text{\rm Aut} (G)$.

%

\underline{Twisted conjugacy problem:} {\it Does there exist an algorithm to
determine is for an arbitrary pair of group  words $U$ and $V$ in generators of
$G$ the equality $\varphi(W) \, U \, = \, V \, W$ holds for some $W \in G$ and
$\varphi \in H$?}

Since the twisted conjugacy problem depends of a group $G$ as well
as a subset $H \subseteq \text{\rm Aut} (G)$, this problem will be
referred to as the conjugacy problem in $G$ {\it relative to} $H$
or as the {\it $H$--twisted conjugacy problem in}~$G$. If $H$
consists of the identity automorphism then we deal with the
classical conjugacy problem. If $H = \{ \varphi \}$ consists of
the unique automorphism then, obviously, a property to be twisted
conjugated is an equivalence relation.

 The following question was posted by G.~Makanin in
``Kourovka notebook'' \cite[Question 10.26(a)]{Kour}: {\it Does
there exist an algorithm to determine is for an arbitrary pair of
words $U$ and $V$ of a free group $G$ and an arbitrary
automorphism $\varphi$ of $G$ the equation $\varphi(X) U = V X$
solvable in $G$?} This question can be regarded as the twisted
conjugacy problem for the case when $G$ is a free group and $H =
\{ \varphi\}$ consists of the unique automorphism.


We remark that the twisted conjugacy problem is also connected
with the question of solving equations in the holomorph $\text{\rm
Hol} (G) = G \ltimes \text{\rm Aut} (G)$. In this context we
recall the following question of  G.~Makanin \cite[Question
10.26(b)]{Kour}: {\it Does there exist an algorithm to determine is
for arbitrary automorphisms $\varphi_1, \ldots, \varphi_n$ of a
free group the equation $w(x_{i_1}^{\varphi_1},\ldots,
x_{i_n}^{\varphi_n})=1$ solvable?}


\section{Basic notations and results}

In the present paper we deal with the case when $G=F_n$ is the
finitely generated free group with basis $\{ x_1, \ldots, x_n\}$.
Elements of $F_n$ are words in the alphabet ${\mathbb X} = \{
x_1^{\pm 1}, \ldots, x_n^{\pm 1}\}$. In our considerations below
we will need to distinguish words which represent the same element
of the group. So, we will usually use small letters to denote
elements of $F_n$ and capital letters to denote words in the above
alphabet. We will write $U = V$ (or $u=v$) if two words (or two
elements of the group) are equal as elements of the group, and $U
\equiv V$ if two words coincide graphically. We use standard
notations $\text{\rm Aut} (F_n)$ and $\text{\rm Inn} (F_n)$ for
the group of automorphisms and the group of inner automorphisms of
$F_n$.

Denote by $|W|$ the length of a word $W$ in the alphabet $\mathbb
X$. A word $W$ is said to be {\it reduced} if it contains no part
$x x^{-1}$, $x \in \mathbb X$. A reduced word $W$ defines a
non-identity element if and only if $|W| \ge 1$. A reduced word
obtained by reducing of an original word will be referred to as
its {\it reduction}. By $|| W ||$ we will denote the length of the
reduction of a word $W$.

We define a linear order ``$<$'' on the set of reduced words in
the alphabet $\mathbb X$. Assume that letters of $\mathbb X$ are
ordered in the following way: $$ x_1 \, < \, x_1^{-1} \, < \, x_2
\, < \, x_2^{-1} \, < \, \ldots \, < \, x_n \, < \, x_n^{-1} .
$$ We write $U < V$ if $|U| < |V|$ or if $|U| = |V|$ and the word
$U$ is less than the word $V$ in respect to the lexicographical
order corresponding to the above defined linear order on $\mathbb
X$.

Recall that all three classical Dehn problems have simple and
elegant solutions in free groups. Since reducing any word is an
algorithmic process, this provides a solution of the word problem.
Further, a reduced word $W=x_{i_1}^{\varepsilon_1}
x_{i_2}^{\varepsilon_2} \cdots x_{i_n}^{\varepsilon_n},$ where
$\varepsilon_i = \pm 1,$ $i=1, \ldots, n$, is said to be {\it
cyclically reduced} if $i_1 \neq i_n$ or if $i_1 = i_n$ then
$\varepsilon_1 \neq - \varepsilon_n$. Clearly, every element of a
free group is conjugated to an element given by a cyclically
reduced word called a {\it cyclic reduction}. This leads to a
solution of the conjugacy problem. Suppose for given words $U$ and
$V$, words $\overline{U}$, $\overline{V}$ are their cyclic
reductions. Then $U$ is conjugated to $V$ if and only if
$\overline{U}$ is a cyclic permutation of $\overline{V}$. Finally,
two finitely generated free groups are isomorphic if and only if
they have the same rank.

In this paper we consider the twisted conjugacy problem in $G=F_n$
relative to $H = \{ \varphi \}$, where $\varphi \in \text{\rm Aut}
(F_n)$ is such that $\varphi^m \in \text{\rm Inn} (F_n)$ for some
non-zero integer $m$. Such automorphism $\varphi$ will be referred
to as a {\it virtually inner automorphism}. In particular, if an
automorphism $\varphi$ is of finite order or inner, then it is
virtually inner. The main result of the paper is the following

\begin{theorem} \label{theorem1}
Suppose $\varphi$ is a virtually inner automorphism of the free
group $F_n$. Then the $\varphi$--twisted conjugacy problem in
$F_n$ is solvable.
\end{theorem}

To prove this statement, we introduce $\varphi$-twisted conjugated
normal form of a word, which can be constructed by a finite number
of steps. We show that two words are $\varphi$-twisted conjugated
if and only if their $\varphi$-twisted conjugated normal form
coincide (see Section~6 for the proof).

\smallskip

Consider the mapping torus (in other words, the ascending HNN
extension)
$$
F_n(\varphi) = \langle x_1, \ldots, x_n, t \ | \ t^{-1} x_i t =
\varphi (x_i), \ i=1, \ldots, n \rangle
$$
that is the semidirect product of $F_n$ and $\langle t \rangle$.
There is an evident relation between the $\varphi$-twisted
conjugacy problem in $F_n$ and the classical conjugacy problem in
$F_n (\varphi)$.

\begin{lemma} \label{lemma1}
The $\varphi$--twisted conjugacy problem in the free group $F_n$
is solvable if and only if there exists an algorithm to determine
is a pair of words $t U$ and $t V$ from $F_n(\varphi)$, where $U,
V \in F_n$, define elements conjugated by some element of $F_n$.
\end{lemma}

\begin{proof}
Obviously, the equality $X^{-1} (t V) X = t U$ holds in
$F_n(\varphi)$ for some $X \in F_n$ if and only if the equality
$\varphi(X^{-1}) V X = U$ holds in $F_n$, i.e. $\varphi (X) U = V
X$.
\end{proof}

It was shown by Bestvina and Feighn in \cite{BF}, that if $\varphi
\in \text{Aut} (F_n)$ have no nontrivial periodic conjugacy
classes (see next section for the definition) then $F_n(\varphi)$
is hyperbolic. So, the conjugacy problem is solvable in it. (We
note, that virtually inner automorphisms that we consider below,
are such that each conjugacy class is periodic.) We will show in
Proposition~\ref{proposition1} that also the problem of
conjugation by an element of $F_n$ is solvable in it.

Theorem~\ref{theorem1}, Lemma~\ref{lemma1}, and
Proposition~\ref{proposition1} imply the partial affirmative
answer on the above mentioned Makanin's question \cite[Question
10.26(a)]{Kour}:

\begin{cor}
Let $\varphi \in \text{Aut} (F_n)$ be a virtually inner or having
no nontrivial periodic conjugacy classes. Then there exists an
algorithm to determine for an arbitrary pair of words $U$ and $V$
the solvability in $F_n$ of the equation $\varphi(X) U = V X$.
\end{cor}

\section{Automorphisms having no nontrivial periodic conjugacy classes}

An automorphism $\varphi$ of the free group $F_n$ is said to have
{\it nontrivial periodic conjugacy class} if there exist integer
$\ell$ and elements $x, y \in F_n$ such that $\varphi^{\ell} (x) =
y^{-1} x y$. It was proven by Bestvina and Feighn \cite{BF} and by
P.~Brinkmann \cite{Br}, that $F_n(\varphi)$ is hyperbolic (in the
Gromov's sense) if and only if automorphism $\varphi$ has no
nontrivial periodic conjugacy classes.

\begin{prop}\label{proposition1}
Suppose $\varphi \in Aut(F_n)$ has no nontrivial periodic
conjugacy classes. Then the $\varphi$-twisted conjugacy problem in
$F_n$ is solvable.
\end{prop}

\begin{proof}
By Lemma~\ref{lemma1}, the $\varphi$--twisted conjugacy of words
$U$ and $V$ from $F_n$ is equivalent to the conjugacy of words
$tU$ and $tV$ from $F_n(\varphi)$ by some element of $F_n$. Since
$\varphi$ has no nontrivial periodic conjugacy classes, by
\cite{BF} the group $F_n(\varphi)$ is hyperbolic and so, the
conjugacy problem in this group is solvable (see, for example,
\cite[Ch.~3]{BH}).

If $tU$ and $tV$ are not conjugated in $F_n(\varphi)$, then they,
in particular, are not conjugated by an element of $F_n$, so $U$
and $V$ are not $\varphi$--twisted conjugated.

Assume that $tU$ and $tV$ are conjugated in $F_n(\varphi)$ by some
element $t^k W$, where $W \in F_n$ and $k$ is integer. If $k=0$,
the statement follows from Lemma~\ref{lemma1}. If $k\neq 0$, we
remark the following. Each element of $F_n(\varphi)$ that
conjugates $tU$ to $tV$ is of the form $c \cdot t^k W$, where $c$
belongs to the centralizer $C_G(tU)$ of $tU$ in $G=F_n(\varphi)$.
It was shown by Gromov (see, for example, \cite[Corollary~3.10,
p.~462]{BH}), that in a hyperbolic group the centralizer of any
nontrivial element is almost cyclic and there exists an algorithm
to find generators of the centralizer $C_G(tU)$. Considering the
canonical homomorphism of the centralizer $C_G(tU)$ onto the group
$\langle t \rangle$, one can check if the element $t^{-k}$ lies in
the image of that homomorphism. If yes, then elements $tU$ and
$tV$ are conjugated by an element of $F_n$. If not, then they are
not conjugated by an element of $F_n$, and so, elements $U$ and
$V$ are not $\varphi$-twisted conjugated.
\end{proof}

We are thankful to Oleg Bogopol'skii for useful discussion of
results from this Section.

\section{Virtually inner automorphisms of free groups}

Let $\varphi \in \text{\rm Aut} (F_n)$ be a virtually inner
automorphism such that $\psi = \varphi^m$ is an inner automorphism
of $F_n$. Without loss of generality we can assume that $m$ is
taken the smallest positive integer having such a property. There
exists a reduced word $\Delta \in F_n$ such that
$$
\varphi^m (f) \, = \, \Delta^{-1} f \Delta
$$
for any $f \in F_n$.


\begin{lemma}\label{lemma2}
(1) Automorphisms $\varphi$ and $\psi$ commute.\\ (2) If $k = m q
+ r$, where $q \in \mathbb Z$ and $0\leq r \leq m-1$, then for any
word $U$ of the alphabet $\mathbb X$ we have
$$
\varphi^k (U) \, = \, \varphi^r (\Delta^{-q} \, U \, \Delta^q).
$$\\
(3) $\varphi (\Delta) = \Delta$.
\end{lemma}

\begin{proof}
Properties (1) and (2) hold obviously, since $\psi$ is a power of $\varphi$.\\
To show (3), remark that by (1) we have
$$
\Delta^{-1} \, \varphi(f) \, \Delta \,  =  \, \varphi (\Delta^{-1}
f \Delta ) \, = \, \varphi(\Delta^{-1}) \, \varphi (f) \,
\varphi(\Delta)
$$
for any $f \in F_n$. Therefore, $\varphi(\Delta) \, \Delta^{-1}$
and $\varphi(f)$ commute for all $f \in F_n$. Since $\varphi$ is
automorphism, we get $\varphi(\Delta) \, \Delta^{-1} = 1$, so
$\varphi(\Delta) = \Delta$.
\end{proof}

A reduced word $V \in F_n$ is said to be {\it $\Delta$-reduced} if
$|V | \leq || \Delta^{-k} V \Delta^k ||$ for all $k \in \mathbb
Z$. Obviously, if $V$ is cyclically reduced, then the length of
any word conjugated to $V$ is not less than the length of $V$, so
$V$ is $\Delta$-reduced.


\begin{lemma}\label{lemma3}
Suppose that $\Delta$ is cyclically reduced. A reduced word $V \in
F_n$ is $\Delta$-{\it reduced} if $|V| \leq ||
\Delta^{-\varepsilon} V \Delta^\varepsilon ||$ for $\varepsilon =
\pm 1$.
\end{lemma}

\begin{proof}
If $V$ is cyclically reduced, by the above observation it is
$\Delta$-reduced. If $V$ is not cyclically reduced, we represent
it in the form $V \equiv U^{-1} W U$, where $U$ and $W$ are
reduced nonempty words in $F_n$ and $W$ is cyclically reduced.

Suppose that $|V| \leq || \Delta^{-1} V \Delta||$ and assume that
there exists integer $k>1$ such that $|| \Delta^{-k} V \Delta^k ||
< | V |$. Then
$$
\Delta^{-k} \, V \, \Delta^k \, \equiv \, \Delta^{-k} \, U^{-1} \,
W \, U \,\Delta^k \, \equiv \, (U \Delta^k)^{-1} \, W \, (U
\Delta^k).
$$
 Since $\Delta$ is cyclically reduced, $|| \Delta^k || =
k \, | \Delta |$ and in the word $U \Delta^k$ only cancellations
of letters from $U$ with letters from $\Delta^k$ can arise. If
such cancellations are possible, there are two possibilities:
either there are cancellations with letters of $W$ or not.

\smallskip

\underline{Case 1.} Suppose that there are cancellations of
letters of $\Delta^k$ with letters of $U$, where, possibly, $U$
will be cancelled wholly (the same for $\Delta^{-k}$ and $U^{-1}$,
respectively), but there no further cancellations with letters of
$W$. Then we can write $U \equiv \Sigma_1 \Sigma_2^{-1}$ and
$\Delta^k \equiv \Sigma_2 \Sigma_3$ for some reduced $\Sigma_1,
\Sigma_2, \Sigma_3 \in F_n$.

If the length of the cancelling part $\Sigma_2$ is bigger than the
length of $\Delta$, then $\Sigma_2 \equiv \Delta \, \Sigma_{2,1}$
and $U \equiv \Sigma_1 \, \Sigma_{2,1}^{-1} \, \Delta^{-1} \equiv
U_1 \, \Delta^{-1}$, where $| U_1 | < |U|$. Then
$$
|| \Delta^{-1} \, V \, \Delta || = || \Delta^{-1} U^{-1} W U
\Delta || = || U_1^{-1} W U_1 || = | W | + 2 | U_1 |
$$
that is less than
$$
|V| = |U^{-1} W U| = |W| + 2 |U|
$$
because $|U_1| < |U|$. Thus we get the contradiction with the
assumption.

If the length of the cancelling part $\Sigma_2$ is less or equal
to the length of  $\Delta$, then $U \equiv \Sigma_1 \Sigma_2^{-1}$
and $\Delta^k \equiv \Sigma_2 \Sigma_3 \Delta^{k-1}$ for some
reduced $\Sigma_1, \Sigma_2, \Sigma_3 \in F_n$. We get
$$
|V| = |W| + 2 |U| = |W| + 2 |\Sigma_1| + 2 |\Sigma_2|,
$$
$$
||\Delta^{-1} V \Delta || = || \Sigma_3^{-1} \Sigma_1^{-1} W
\Sigma_1 \Sigma_3 || = |W| + 2 |\Sigma_1| + 2|\Sigma_3|,
$$
and
$$ || \Delta^{-k} V \Delta^k || = || \Delta^{-(k-1)}
(\Delta^{-1} V \Delta) \Delta^{k-1} ||$$
$$
= || \Delta^{-1} V \Delta || + 2 (k-1) |\Delta |  = |W| + 2
|\Sigma_1| + 2|\Sigma_3| + 2 (k-1) |\Delta|.
$$
Since $k > 1$ and $|\Sigma_2| \leq |\Delta|$, we see that $|V|
\leq ||\Delta^{-k} V \Delta^k ||$ that gives the contradiction
with the assumption.

\smallskip

\underline{Case 2.} Suppose that letters of $\Delta^{-k}$ and
$\Delta^k$ are cancelling with letters of $V \equiv U^{-1} W U$ is
such a way that words $U$ and $U^{-1}$ will be cancelled wholly
and there are further cancellations with letters of $W$. Since
$|U^{-1} W U| \leq ||\Delta^{-1} U^{-1} W U \Delta||$, and
$U^{-1}$ and $U$ are cancelling wholly, we have $U \equiv
\Delta_1^{-1}$ and $\Delta \equiv \Delta_1 \, \Delta_2$ with
$|\Delta_2| \geq |\Delta_1| = |U|$. Hence
\begin{eqnarray*}
\Delta^{-k} \, V \, \Delta^k & \equiv &  \Delta^{-(k-1)} \,
\Delta_2^{-1} \, \Delta_1^{-1} \, (\Delta_1 \, W \, \Delta_1^{-1}
) \, \Delta_1 \, \Delta_2 \, \Delta^{k-1} \cr & = &
\Delta^{-(k-1)} \, \Delta_2^{-1} \, W \, \Delta_2 \, \Delta^{k-1}.
\end{eqnarray*}
The further cancellations are possible either in $\Delta_2^{-1} \,
W$ or in $W \Delta_2$, but not in both, since $W$ is cyclically
reduced.

If $\Delta^{-1}_2$ is cancelling wholly in $\Delta_2^{-1} W
\Delta_2$, then
$$
|| \Delta^{-1} V \Delta ||  = || \Delta^{-1}_2 W \Delta_2 ||  =
|W| < |V|
$$
and we have a contradiction with the assumption $||\Delta^{-1} V
\Delta|| \geq |V|$.

If $\Delta^{-1}_2$ is not cancelling wholly in the product
$\Delta^{-1}_2 W \Delta_2$, then $\Delta_2 \equiv \Delta_{21}
\Delta_{22}$ and $W \equiv \Delta^{-1}_{21} W_1$, and
\begin{eqnarray*}
|| \Delta^{-k} V \Delta^k || & = & || \Delta^{-(k-1)}
\Delta_2^{-1} W \Delta_2 \Delta^{k-1} || \cr & = & ||
\Delta^{-(k-1)} \Delta^{-1}_{22} W_1 \Delta_{21} \Delta_{22}
\Delta^{k-1} || \cr &
\geq & |W| + 2 (k-1) \, |\Delta| + 2 |\Delta_{22}| \\
& > & |W| + 2 \, |\Delta| \, > \, |W| + 2 |U| = |V|,
\end{eqnarray*}
and we have a contradiction with the assumption $||\Delta^{-k} V
\Delta^k|| < |V|$. By similar arguments, the case when $\Delta_2$
is not cancelling wholly in the product $W \Delta_2$, is also
impossible.

\smallskip

If there exists integer $k<-1$ with the same property, similar
considerations implies the statement.
\end{proof}

If $\Delta$ is cyclically reduced, Lemma~\ref{lemma3} gives the
finite algorithm to find for a given reduced word $V$ a
$\Delta$-reduced word $V_{\Delta}$ conjugated to $V$ by some power
of $\Delta$. Indeed, it is enough to repeat conjugations of $V$ by
$\Delta^{\varepsilon}$, $\varepsilon = \pm 1$, few times. If the
length of the obtained word is less than the length of the
previous word, we will conjugate again. If not, then the obtained
word is a $\Delta$-reduced word $V_{\Delta}$ conjugated to $V$.
Such a construction of a $\Delta$-reduced word $V_{\Delta}$
conjugated to $V$ by some power of $\Delta$ will be referred to as
a $\Delta$-{\it reduction}.

We remark that if $\Delta$ is not cyclically reduced, then the
analog of Lemma~\ref{lemma3} does not hold. It is clear from the
following example.

\smallskip

{\bf Example.} Let $U, W, \Sigma \in F_n$ be nonempty reduced
words such that for an integer $|k|>1$ words $\Delta \equiv U^{-1}
W^{-1} \Sigma W U$ and $V \equiv U^{-1} W^{-1} \Sigma^k W U^2$ are
reduced. If $k>1$, we get
\begin{eqnarray*}
\Delta^{-1} V \Delta & = & U^{-1} W^{-1} \Sigma^{k-1} W U W^{-1}
\Sigma W U, \cr \Delta^{-2} V \Delta^2 & = & U^{-1} W^{-1}
\Sigma^{k-2} W U W^{-1} \Sigma^2 W U, \cr & \cdots & \cr
\Delta^{-(k-1)} V \Delta^{k-1} & = & U^{-1} W^{-1} \Sigma W U
W^{-1} \Sigma^{k-1} W U, \cr \Delta^{-k} V \Delta^k & = & W^{-1}
\Sigma^k W U. \end{eqnarray*}
It is easy to see that
$$
|V| < ||\Delta^{-1} V \Delta|| = ||\Delta^{-2} V \Delta^2 || =
\ldots = ||\Delta^{-(k-1)} V \Delta^{k-1}||,
$$
but $||\Delta^{-k} V \Delta^k || < |V|$.

If $k < -1$, similar example can be obtained conjugating $V$ by
negative powers of $\Delta$.



\smallskip

\begin{lemma}\label{lemma4}
Suppose $\Delta$ is not cyclically reduced. Let $V$ be a reduced
word such that $|V| \leq ||\Delta^{-\varepsilon} V
\Delta^{\varepsilon} ||$, $\varepsilon = \pm 1$. Then one of the following cases
holds:\\
(1) $V$ is $\Delta$-reduced.  \\
(2) $\Delta \equiv U_1^{-1} U_2^{-1} \Delta_{11}^{-1} \Delta_2
\Delta_{11} U_2 U_1$ and $V \equiv U_1^{-1} U_2^{-1} W U_2 U_1$,
for some reduced $\Delta_{11}, \Delta_{2}, U_1, U_2$ and
cyclically reduced $W$ for which there exist integers $k$,
$|k|>1$, and
$m \geq 0$ and reduced $\Phi$, $W_0$ such that either\\
a) $U_2^{-1} \Delta_{11}^{-1} \Delta_2^{-k} \Delta_{11} \equiv
\Phi W^{-m}$ and $W \equiv \Phi^{-1} W_0$,\\
or\\
b) $\Delta_{11}^{-1} \Delta_2^{k} \Delta_{11} U_2 \equiv W^{-m}
\Phi$ and $W \equiv W_0 \Phi^{-1}$.\\
In addition, $\Phi$ does not end by $W^{-1}$, and $W_0$ is either
nonempty or $W_0 \equiv 1$ with $U_1^{-1} \Phi^{-1} U_1$ be
reduced.

In case (2) $V_{\Delta} = \Delta^{-k} V \Delta^k$ is
$\Delta$-reduced word conjugated to $V$. Moreover, case (2)
describe all possible cases when $|V|\leq ||\Delta^{-\varepsilon}
V \Delta^{\varepsilon}||$, $\varepsilon = \pm 1$, but
$||\Delta^{-k} V \Delta^k|| < |V|$ for some $|k|>1$.
\end{lemma}

\begin{proof}
Suppose that $|V| \leq ||\Delta^{-1} V \Delta||$ and assume that
there exists integer $k>1$ such that $||\Delta^{-k} V \Delta^k|| <
|V|$.

We represent the reduced word $\Delta$ in the form $\Delta \equiv
\Delta_1^{-1} \Delta_2 \Delta_1$, where $\Delta_2$  is cyclically
reduced and $\Delta_1$ is reduced. Then $\Delta^k  =
(\Delta_1^{-1} \Delta_2 \Delta_1)^k = \Delta_1^{-1} \Delta_2^k
\Delta_1$, where the last word is reduced and $||\Delta^k|| = 2
|\Delta_1| + k |\Delta_2|$.

If $V$ is cyclically reduced, then, obviously, it is
$\Delta$-reduced.

If $V$ is not cyclically reduced, we write it in the form $V\equiv
U^{-1} W U$, where $W$ is cyclically reduced and $U$ is reduced.
Then
$$
\Delta^{-k} V \Delta^k = \Delta^{-k} U^{-1} W U \Delta^k = (U
\Delta_1^{-1} \Delta_2^k \Delta_1)^{-1} W (U \Delta_1^{-1}
\Delta_2^k \Delta_1) .
$$
Since $|| \Delta^k V \Delta^{-k}|| < |V|$, there are cancellations
in this word. There are two possibilities: either there are
cancellations with letters of $W$ or not.

\smallskip

\underline{Case 1.} Suppose that there are cancellations of
letters of $\Delta^k$ with letters of $U$, where, possibly, $U$
will be cancelled wholly (the same for $\Delta^{-k}$ and $U^{-1}$,
respectively), but there no further cancellations with letters of
$W$. Then we can write $U \equiv \Sigma_1 \Sigma_2^{-1}$ and
$\Delta_1^{-1} \Delta_2^k \Delta_1 \equiv \Sigma_2 \Sigma_3$ for
some $\Sigma_1, \Sigma_2, \Sigma_3 \in F_n$.

Since
$$
|| \Delta^{-k} V \Delta^k ||  = | \Sigma_3^{-1} \Sigma_1^{-1} W
\Sigma_1 \Sigma_3 | = |W| + 2 | \Sigma_1 | + 2 | \Sigma_3 |
$$
and
$$
| V | = | \Sigma_2 \Sigma_1^{-1} W \Sigma_1 \Sigma_2^{-1} | = |W|
+ 2 | \Sigma_1 | + 2 |\Sigma_2 |,
$$
we get $|\Sigma_2 | > |\Sigma_3|$.

{}From the representation $\Delta_1^{-1} \Delta_2^k \Delta_1
\equiv \Sigma_2 \Sigma_3$ we have $\Sigma_2 \equiv \Delta_1^{-1}
\Delta_2^{\ell} \Delta_{21}$ and $\Sigma_3 \equiv \Delta_{22}
\Delta_2^m \Delta_1$, where $\Delta_2 \equiv \Delta_{21}
\Delta_{22}$ (one of words $\Delta_{21}$ and $\Delta_{22}$ can be
empty) and $k = \ell +m + 1$.

Consider
\begin{eqnarray*}
|| \Delta^{-1} V \Delta || & = & || \Delta_1^{-1} \Delta_2^{-1}
\Delta_1 U^{-1} W U \Delta_1^{-1} \Delta_2 \Delta_1 || \cr & = &
|| \Delta_1^{-1} \Delta_2^{-1} \Delta_1 \Sigma_2 \Sigma_1^{-1} W
\Sigma_1 \Sigma_2^{-1} \Delta_1^{-1} \Delta_2 \Delta_1 || \cr & =
& || \Delta_1^{-1} \Delta_2^{-1} \Delta_1 ( \Delta_1^{-1}
\Delta_2^{\ell} \Delta_{21} ) \Sigma_1^{-1} W \Sigma_1 (
\Delta_{21}^{-1} \Delta_2^{-\ell} \Delta_1 ) \Delta_1^{-1}
\Delta_2 \Delta_1 || \cr & = & || \Delta_1^{-1} \Delta_2^{\ell -1}
\Delta_{21} \Sigma_1^{-1} W \Sigma_1 \Delta_{21}^{-1}
\Delta_2^{-\ell + 1} \Delta_1 || \cr & = & |W| + 2 \, |\Sigma_1| +
2 \, |\Delta_1^{-1} \Delta_2^{\ell-1} \Delta_{21} |
\end{eqnarray*}
and
$$
| V | = | \Sigma_2 \Sigma_1^{-1} W \Sigma_1 \Sigma_2^{-1} | = | W
| + 2 | \Sigma_1 | + 2 | \Sigma_2 | .
$$
Since $\Sigma_2 \equiv \Delta_1^{-1} \Delta_2^{\ell} \Delta_{21}$,
we have $|\Sigma_2| > |\Delta_1^{-1} \Delta_2^{\ell-1} \Delta_{21}
|$. Therefore, $||\Delta^{-1} V \Delta || < | V |$, that gives the
contradiction with the assumption.

Therefore, for any $k \ge 1$ the inequality $|V| \leq ||
\Delta^{-k} V \Delta^k ||$ holds.

Taking $\Delta V \Delta^{-1}$, similar considerations show that
the inequality holds also for $k \leq -1$. Therefore, $V$ is
$\Delta$-reduced word.

\smallskip

\underline{Case 2.} Suppose that letters of $\Delta^{-k}$ and
$\Delta^k$, where $\Delta^k \equiv \Delta_1^{-1} \Delta_2^k
\Delta_1$ are cancelling with letters of $V \equiv U^{-1} W U$ is
such a way that words $U$ and $U^{-1}$ will be cancelled wholly
and there are further cancellations with letters of $W$, starting
either from the initial part of $W$ or from the final part of $W$,
but not from the both, since $W$ is cyclically reduced.

\smallskip

\underline{Case 2(i).} Suppose that $|U| \geq |\Delta_1|$ and
cancellations in $W$ starts from the initial part.

Recall that $|V| \leq ||\Delta^{-1} V \Delta||$, where $|V| = |W|
+ 2 |U|$ and $||\Delta^{-1} V \Delta|| = || \Delta_1^{-1}
\Delta_2^{-1} \Delta_1 U^{-1} W U \Delta_1^{-1} \Delta_2 \Delta_1
||$. Hence $|U| \leq |\Delta_1| + \frac{1}{2} |\Delta_2|$ and we
can represent $U^{-1} \equiv \Delta_1^{-1} \Delta_{21}$, where
$\Delta_2 \equiv \Delta_{21} \Delta_{22}$ and $|\Delta_{21}| \leq
|\Delta_{22} |$. Then
\begin{eqnarray*}
\Delta^{-k} V \Delta^k & = & \Delta_1^{-1} (\Delta_{21}
\Delta_{22})^{-k} \Delta_1 (\Delta_1^{-1} \Delta_{21}) W
(\Delta_{21}^{-1} \Delta_1) \Delta_1^{-1} (\Delta_{21}
\Delta_{22})^k \Delta_1 \\ & = & \Delta_1^{-1} (\Delta_{21}
\Delta_{22})^{-(k-1)} \Delta_{22}^{-1} W \Delta_{22} (\Delta_{21}
\Delta_{22})^{k-1} \Delta_1 .
\end{eqnarray*}
Consider
$$
|V| = |U^{-1} W U | = |W| + 2 |U| = |W| + 2 |\Delta_1| + 2 |
\Delta_{21}|
$$
and
$$
\Delta^{-1} V \Delta = \Delta_1^{-1} \Delta_2^{-1} \Delta_1 U^{-1}
W U \Delta_1^{-1} \Delta_2 \Delta_1  = \Delta_1^{-1}
\Delta_{22}^{-1} W \Delta_{22} \Delta_1 .
$$
Since $|V| \; \leq || \Delta^{-1} V \Delta ||$, in the product
$\Delta_{22}^{-1} W$ the word $\Delta_{22}^{-1}$ can not be
cancelled wholly.

Indeed, if $\Delta_{22}^{-1}$ is cancelling wholly, then after
cancellations in $\Delta_{22}^{-1} W \Delta_{22}$ we will get a
cyclically reduced word which is obtained by a cyclic shift  of
$W$. Therefore $||\Delta_{22}^{-1} W \Delta_{22}|| = |W|$ and
$$
|| \Delta^{-1} V \Delta || \leq  |W| + 2 | \Delta_1| < |V|,
$$
and we will get the contradiction.

Therefore, we can represent $\Delta_{22}^{-1} \equiv \Sigma_0
\Sigma_1^{-1}$, where $\Sigma_1^{-1}$ is the cancelling part and
$\Sigma_0$ is non-empty. After cancellations in $\Sigma_1^{-1} W
\Sigma_1$ we will get a cyclically reduced word which is obtained
by a cyclic shift of $W$. Therefore, $||\Sigma_1^{-1} W \Sigma_1
|| = |W|$ and
$$
|| \Delta^{-1} V \Delta || = || \Delta_1^{-1} \Sigma_0
(\Sigma_1^{-1} W \Sigma_1) \Sigma_0^{-1} \Delta_1 || = |W| + 2
|\Delta_1| + 2 |\Sigma_0| .
$$
Comparing with
\begin{eqnarray*} || \Delta^{-k} V \Delta^k || & = & ||
\Delta_1^{-1} (\Delta_{21} \Delta_{22})^{-(k-1)} \Delta_{22}^{-1}
W \Delta_{22} (\Delta_{21} \Delta_{22})^{k-1} \Delta_1 || \\
& = & || \Delta_1^{-1} (\Delta_{21} \Delta_{22})^{-(k-1)} \Sigma_0
(\Sigma_1^{-1} W \Sigma_1) \Sigma_0^{-1}
(\Delta_{21} \Delta_{22})^{k-1} \Delta_1 || \\
& = & | W | + 2 | \Delta_1 | + 2 | \Sigma_0 | + 2(k-1) |
\Delta_{21} \Delta_{22} | ,
\end{eqnarray*}
we see that $||\Delta^{-1} V \Delta|| < ||\Delta^{-k} V
\Delta^k||$ for all $k > 1$. By similar considerations the
inequality holds also for all $k < -1$. Therefore, $V$ is
$\Delta$-reducible.

\smallskip

\underline{Case 2(ii).} Suppose that $|U| \geq |\Delta_1|$ and
cancellations in $W$ starts from the final part. The arguments
analogous to the arguments from the Case~2(i), shows that in this
case $V$ $\Delta$-reducible.

\smallskip

\underline{Case 2(iii).} Suppose that $|U| < |\Delta_1|$ and
cancellations in $W$ starts from the initial part. Then we can
write $\Delta_1 \equiv \Delta_{11} U$, where $\Delta_{11} \neq 1$.
Therefore,
\begin{eqnarray*}
\Delta^{-k} V \Delta^k & \equiv & \Delta_1^{-1} \Delta_2^{-k}
\Delta_1 U^{-1} W U \Delta_1^{-1} \Delta_2^k \Delta_1 \\
& = & U^{-1} \Delta_{11}^{-1} \Delta_2^{-k} \Delta_{11} U U^{-1} W
U U^{-1} \Delta_{11}^{-1} \Delta_2^{k} \Delta_{11} U \\
& = & U^{-1} \Delta_{11}^{-1} \Delta_2^{-k} \Delta_{11} W
\Delta_{11}^{-1} \Delta_2^k \Delta_{11} U .
\end{eqnarray*}
Since $||\Delta^{-k} V \Delta^k|| < |V|$, where $|V| = |W| + 2
|U|$, is it necessary that in the product $U^{-1} \Delta_{11}^{-1}
\Delta_2^{-k} \Delta_{11} W$ the word $\Delta_{11}^{-1}
\Delta_2^{-k} \Delta_{11}$ will be cancelled wholly and, moreover,
there will be some cancellations of letters of $U^{-1}$ with
letters of $W$. Denote by $U_2^{-1}$ the cancelling part of the
word $U^{-1}$, Then we can write $U \equiv U_2 U_1$, where $|U| =
|U_1| + |U_2|$, and $W \equiv \Delta_{11}^{-1} \Delta_2^k
\Delta_{11} U_2 W_0$, for some $W_0$. Since $W$ is cyclically
reduced, the word $W_0 \Delta_{11}^{-1} \Delta_2^k \Delta_{11}
U_2$ is reduced and $|| W_0 \Delta_{11}^{-1} \Delta_2^k
\Delta_{11} U_2|| = |W|$. Hence
\begin{eqnarray*} || \Delta^{-k} V \Delta^k || & = & || U_1^{-1} U_2^{-1}
\Delta_{11}^{-1} \Delta_2^{-k} \Delta_{11} W \Delta_{11}^{-1}
\Delta_2^k \Delta_{11} U_2 U_1 || \\
& = & || U_1^{-1} W_0 \Delta_{11}^{-1} \Delta_2^k \Delta_{11} U_2
U_1 || \\
& = & |W| + 2 |U_1| \, < \, |V|,
\end{eqnarray*}
because $|U_1| < |U|$.

Remark that with such $V$ and $\Delta$ the inequality
$||\Delta^{-1} V \Delta|| > |V|$ also holds. Indeed,
\begin{eqnarray*}
\Delta^{-1} V \Delta & \equiv & U^{-1} \Delta_{11}^{-1}
\Delta_2^{-1} \Delta_{11} W \Delta_{11}^{-1} \Delta_2 \Delta_{11}
U \\
& \equiv & U_1^{-1} U_2^{-1} \Delta_{11}^{-1} \Delta_2^{-1}
\Delta_{11} ( \Delta_{11}^{-1} \Delta_2^k \Delta_{11} U_2 W_0 )
\Delta_{11}^{-1} \Delta_2 \Delta_{11} U_2 U_1 \\
& = & U_1^{-1} U_2^{-1} \Delta_{11}^{-1} \Delta_2^{(k-1)}
\Delta_{11} U_2 W_0 \Delta_{11}^{-1} \Delta_2 \Delta_{11} U_2 U_1.
\end{eqnarray*}
Since the obtained word is reduced, we get $||\Delta^{-1} V \Delta
|| = |W| + 2 |U| + 2 |\Delta_{11}|$, but $|V| = |W| + 2 |U|$.

Also, it is easy to see that $||\Delta V \Delta^{-1}||
> |V|$.

Thus we get that if
$$
\Delta \equiv U_1^{-1} U_2^{-1} \Delta_{11}^{-1} \Delta_2
\Delta_{11} U_2 U_1
$$
and
$$
V \equiv
U_1^{-1} U_2^{-1} \Delta_{11}^{-1} \Delta_2^k \Delta_{11} U_2 W_0
U_2 U_1
$$ for some $\Delta_{11}, \Delta_2, U_1, U_2, W_0 \in F_n$,
and some integer $k>1$, then $||\Delta^{-\varepsilon} V
\Delta^{\varepsilon} || > |V|$, but $||\Delta^{-k} V \Delta^k|| <
|V|$. Note that we have $W = \Delta_{11}^{-1} \Delta_2^k
\Delta_{11} U_2 W_0$, that gives the case (2a) of the statement
for the case $m=0$.

Above we considered the case when in the product $U_2^{-1}
\Delta_{11}^{-1} \Delta_2^{-k} \Delta_{11} W$ the word $U_2^{-1}
\Delta_{11}^{-1} \Delta_2^{-k} \Delta_{11}$ was cancelling wholly
with initial part of $W$ of with the whole $W$, but without
further cancellations. Now, let us consider the case when $W$ is
cancelling wholly and after that there are further cancellations
in the product
$$ \Delta^{-k} V \Delta^k = U_1^{-1} U_2^{-1}
\Delta_{11}^{-1} \Delta_2^{-k} \Delta_{11} W \Delta_{11}^{-1}
\Delta_2^{k} \Delta_{11} U_2 U_1 .
$$
Remark, that it is possible only if the word $U_2^{-1}
\Delta_{11}^{-1} \Delta_2^{-k} \Delta_{11}$ is a product of some
initial subword $\Phi$ and of an element from the centralizer of
$W$, i.e.
$$
U_2^{-1} \Delta_{11}^{-1} \Delta_2^{-k} \Delta_{11} \equiv \Phi
W^{-m}, \quad m \in \mathbb{N}
$$
In this case we assume that $m$ is maximal integer with such
property, i.e. $\Phi$ does not contain word $W^{-1}$ as a final
part. Then
$$
\Delta^{-k} V \Delta^k = U_1^{-1} \Phi W^{-m} W W^m \Phi^{-1} U_1
= U_1^{-1} \Phi W \Phi^{-1} U_1
$$
and after that $\Phi$ is cancelling with an initial subword of
$W$, i.e. $W \equiv \Phi^{-1} W_0$, where $W_0$ can, possibly, be
empty. If $W_0$ is empty, then, obviously, $U_1^{-1} \Phi^{-1}
U_1$ must be reduced.

\smallskip

\underline{Case 2(iv).} Suppose that $|U| < |\Delta_1|$ and
cancellations in $W$ starts from the final part. By the same
arguments as in the Case~2(iii), we get that words $\Delta$ and
$V$ are of the form, described in case (2b).

\end{proof}

If $\Delta$ is not cyclically reduced, Lemma~\ref{lemma4} gives
the finite algorithm to find for a given reduced word $V$ a
$\Delta$-reduced word $V_{\Delta}$ conjugated to $V$ by some power
of $\Delta$. If $\Delta$ and $V$ are of the form represented in
case (2) of Lemma~\ref{lemma4}, then we define $V_{\Delta} =
\Delta^{-k} V \Delta^k$. In this case we say that $V_{\Delta}$ is
obtained from $V$ by $\Delta$-reducing. If $\Delta$ and $V$ are
others, then we follow the same steps as described after
Lemma~\ref{lemma3}.

\section{$\Delta$-reduced words corresponding to the same word}

In general, $V_{\Delta}$ is not uniquely determined by $V$. The
following statement describes different $\Delta$-reduced words
corresponding to the same reduced word.

\begin{prop} \label{proposition2}
(1) Let $V$ be reduced and $| V |  = || \Delta^{-1} V \Delta||$.
Then $$\Delta \equiv \Delta_1 \Delta_{31} \Delta_{21} \Delta_{31}
\Delta_{32},$$ where $|\Delta_1| = |\Delta_{32}|$, with either
$$V \equiv \Delta_1 V_0 \Delta_{21}^{-1} \Delta_{31}^{-1}
\Delta_1^{-1} \quad \text{or} \quad V \equiv \Delta_1 \Delta_{31}
\Delta_{21} V_0 \Delta_1^{-1},$$ where
$\Delta_{31} \equiv 1$ if  $V_0 \neq 1$, for some $\Delta_1,
\Delta_{21}, \Delta_{31}, \Delta_{32}, V_0 \in F_n$.\\
(2) Let $V$ be reduced and $| V | = || \Delta V \Delta^{-1}||$.
Then $$\Delta \equiv \Delta_{11} \Delta_{12} \Delta_{21}
\Delta_{12} \Delta_{3},$$ where $|\Delta_{11}| = |\Delta_{3}|$,
with either $$V \equiv \Delta_3^{-1} V_0 \Delta_{21} \Delta_{12}
\Delta_3 \quad \text{or} \quad  V \equiv \Delta_3^{-1}
\Delta_{12}^{-1} \Delta_{21}^{-1} V_0 \Delta_3,$$ where
$\Delta_{12} \equiv 1$ if $V_0 \neq 1$, for some $\Delta_{11},
\Delta_{12}, \Delta_{21}, \Delta_{3}, V_0 \in F_n$.
\end{prop}

\begin{proof} Let us prove the statement (1).
Represent a $\Delta$-reduced word $V$ in the form $V \equiv U^{-1}
W U$, where $W$ is cyclically reduced and $U$ is reduced (we
assume that $U$ is empty if $V$ is cyclically reduced itself, i.e.
$V\equiv W$). Since $| V | = || \Delta V \Delta^{-1}||$, the word
$\Delta^{-1} V \Delta \equiv \Delta^{-1} U^{-1} W U \Delta$ admits
cancellations. There are two possibilities: either there are
cancellations with letters of $W$ or not.

\smallskip

\underline{Case 1:} Suppose that there are cancellations of
letters of $\Delta$ with letters of $U$, where, possibly, $U$ will
be cancelled wholly (the same for $\Delta^{-1}$ and $U^{-1}$,
respectively), but there no further cancellations with letters of
$W$. Then $U \equiv U_1 \Delta_1^{-1}$ and $\Delta \equiv \Delta_1
\Delta_3$, where the reduced word $\Delta_1$ is non-empty and
reduced words $U_1$ and $\Delta_3$ are, possibly, empty.
Therefore,
$$
\Delta^{-1}  V  \Delta  \equiv \Delta_3^{-1} \Delta_1^{-1} (
\Delta_1 U_1^{-1} W U_1 \Delta_1^{-1} ) \Delta_1 \Delta_3 =
\Delta_3^{-1} U_1^{-1} W U_1 \Delta_3
$$
where the final word is reduced. So,
$$
|| \Delta^{-1} V \Delta || = | W | + 2 |U_1| + 2 |\Delta_3| .
$$
Since
$$
| V | = | \Delta_1 U_1^{-1} W U_1 \Delta_1^{-1} | =  |W| + 2 |U_1|
+ 2 |\Delta_1|,
$$
we get $|\Delta_1| = |\Delta_3|$. Taking $V_0 \equiv U_1^{-1} W
U_1$, $\Delta_{32} \equiv \Delta_3$, $\Delta_{21} \equiv 1$ and
$\Delta_{31} \equiv 1$ we get $\Delta$ and $V$ of the same form as
in the statement.

\smallskip

\underline{Case 2:} Suppose that words $U$ and $U^{-1}$ will be
cancelled wholly and there are further cancellations with letters
of $W$, starting either from the initial part of $W$ or from the
final part of $W$, but not from both, because $W$ is cyclically
reduced.

\smallskip

\underline{Case 2(i):} Suppose that there are cancellations
starting from the final part of $W$ such that, possibly, $W$
cancelled wholly, but there are no further cancellations. Then
$\Delta \equiv \Delta_1 \Delta_2 \Delta_3$, $U \equiv
\Delta_1^{-1}$ and $W \equiv V_0 \Delta_2^{-1}$ for some reduced
$\Delta_1, \Delta_2, \Delta_3, V_0 \in F_n$. Therefore, $V \equiv
\Delta_1 V_0 \Delta_2^{-1} \Delta_1^{-1}$ and $\Delta^{-1} V
\Delta = \Delta_3^{-1} \Delta_2^{-1} V_0 \Delta_3$. Since
$$
|V|  = |V_0| + |\Delta_2| + 2 |\Delta_1|
$$
and
$$
|| \Delta^{-1} V \Delta || = |V_0| + |\Delta_2| + 2 |\Delta_3| ,
$$
we get $|\Delta_1| = | \Delta_3|$. Taking $\Delta_{21} \equiv
\Delta_2$, $\Delta_{32} \equiv \Delta_3$ and $\Delta_{31} \equiv
1$ we get $\Delta$ and $V$ of the same form as in the statement
(with $V_0 \equiv 1$ if $W$ cancelled wholly).

\smallskip

\underline{Case 2(ii):} Suppose that $W$ is cancelled wholly
starting from the final part and after that some cancellations are
still possible. Therefore, $V_0 \equiv 1$, $\Delta^{-1} V \Delta
\equiv \Delta_3^{-1} \Delta_2^{-1} \Delta_3$ and some
cancellations of letters of $\Delta_2^{-1}$ and $\Delta_3$ are
possible. Thus, denoting by $\Delta_{32}$ the part (possibly,
empty) of $\Delta_3$ which cannot be cancelled, we can write
$\Delta_3 \equiv \Delta_{31} \Delta_{32}$ and $\Delta^{-1} V
\Delta = \Delta_{32}^{-1} \Delta_{31}^{-1} \Delta_2^{-1}
\Delta_{31} \Delta_{32}$, where $\Delta_{31}^{-1} \Delta_2^{-1}
\Delta_{31}$, after possible cancellations, is cyclically reduced.
So, elements $\Delta_2^{-1}$ and $\Delta_{31}^{-1} \Delta_2^{-1}
\Delta_{31}$ are conjugated in the free group and cyclically
reduced. Therefore, the reduction of the word $\Delta_{31}^{-1}
\Delta_2^{-1} \Delta_{31}$ is a cyclic permutation of the word
$\Delta_2^{-1}$, so $\Delta_2^{-1} = \Delta_{21}^{-1}
\Delta_{31}^{-1}$ and $\Delta_{31}^{-1} \Delta_2^{-1}
\Delta_{31}^{-1} = \Delta_{31}^{-1} \Delta_{21}^{-1}
\Delta_{31}^{-1} \Delta_{31}  = \Delta_{31}^{-1}
\Delta_{21}^{-1}$. Hence
$$
\Delta \equiv \Delta_1 \Delta_2 \Delta_3 \equiv \Delta_1
\Delta_{31} \Delta_{21} \Delta_{31} \Delta_{32},
$$
$$
V \equiv \Delta_1 \Delta_2^{-1} \Delta_1^{-1} \equiv \Delta_1
\Delta_{21}^{-1} \Delta_{31}^{-1} \Delta_1^{-1}
$$
and
$$
\Delta^{-1} V \Delta = \Delta_3^{-1} \Delta_2^{-1} \Delta_3 =
\Delta_{32}^{-1} \Delta_{31}^{-1} \Delta_{21}^{-1}
\Delta_{31}^{-1} \Delta_{31} \Delta_{32} = \Delta_{32}^{-1}
\Delta_{31}^{-1} \Delta_{21}^{-1} \Delta_{32} .
$$
Since
$$
|V| = 2 |\Delta_1| + |\Delta_{21}| + |\Delta_{31}|
$$
and
$$
|| \Delta^{-1} V \Delta || = 2 |\Delta_{32}| + |\Delta_{21}| +
|\Delta_{31}|,
$$
we get $|\Delta_1| = |\Delta_{32}|$. Taking $V_0 \equiv 1$ we get
that $\Delta$ and $V$ are of the same form as in the statement.

\smallskip

\underline{Case 2(iii):} Suppose that there are cancellations
starting from the initial part of $W$ such that, possibly, $W$
cancelled wholly, but there no further cancellations. Then $\Delta
\equiv \Delta_1 \Delta_2 \Delta_3$, $U \equiv \Delta_1^{-1}$ and
$W \equiv \Delta_2 V_0$ for some reduced words $\Delta_1,
\Delta_2, \Delta_3, V_0 \in F_n$. Therefore, $V \equiv \Delta_1
\Delta_2 V_0 \Delta_1^{-1}$ and $\Delta^{-1} V \Delta =
\Delta_3^{-1} V_0 \Delta_2 \Delta_3$. Since
$$
|V| = |V_0| + |\Delta_2| + 2 |\Delta_1|
$$
and
$$
|| \Delta^{-1} V \Delta || = |V_0| + |\Delta_2| + 2 |\Delta_3|,
$$
we get $|\Delta_1|  = |\Delta_3|$. Taking $\Delta_{21} \equiv
\Delta_2$, $\Delta_{32} \equiv \Delta_3$ and $\Delta_{31} \equiv
1$ we get that $\Delta$ and $V$ are of the form as in the
statement.

\smallskip

\underline{Case 2(iv):} Suppose that $W$ is cancelled wholly
(starting from the initial part) and after that some cancellations
are still possible. Therefore, $V_0 \equiv 1$, $\Delta^{-1} V
\Delta = \Delta_3^{-1} \Delta_2 \Delta_3$ and some cancellations
of letters of $\Delta_3^{-1}$ and $\Delta_2$ are possible.
Denoting by $\Delta_{32}^{-1}$ the part (possibly, empty) of
$\Delta_3^{-1}$ which can not be cancelled, we can write
$\Delta_3^{-1} \equiv \Delta_{32}^{-1} \Delta_{31}^{-1}$, so
$\Delta_3 \equiv \Delta_{31} \Delta_{32}$. Therefore we get
$\Delta^{-1} V \Delta = \Delta_3^{-1} \Delta_2 \Delta_3 =
\Delta_{32}^{-1} \Delta_{31}^{-1} \Delta_2 \Delta_{31}
\Delta_{32}$, where $\Delta_{31}^{-1} \Delta_2 \Delta_{31}$, after
possible cancellations, is cyclically reduced. Elements $\Delta_2$
and $\Delta_{31}^{-1} \Delta_2 \Delta_{31}$ are conjugated in the
free group and cyclically reduced. Therefore, the reduction of the
word $\Delta_{31}^{-1} \Delta_2 \Delta_{31}$ is a cyclic
permutation of the word $\Delta_2^{-1}$, so $\Delta_2 =
\Delta_{31} \Delta_{21}$ and $\Delta_{31}^{-1} \Delta_2
\Delta_{31} = \Delta_{31}^{-1} \Delta_{31} \Delta_{21} \Delta_{31}
= \Delta_{21} \Delta_{31}$. Hence
$$
\Delta \equiv \Delta_1 \Delta_2 \Delta_3 \equiv \Delta_1
\Delta_{31} \Delta_{21} \Delta_{31} \Delta_{23},
$$
$$
V \equiv \Delta_1 \Delta_2 \Delta_1^{-1} \equiv \Delta_1
\Delta_{31} \Delta_{21} \Delta_1^{-1},
$$
and
$$
\Delta^{-1} V \Delta = \Delta_3^{-1} \Delta_2 \Delta_3 =
\Delta_{32}^{-1} \Delta_{31}^{-1} \Delta_{31} \Delta_{21}
\Delta_{31} \Delta_{32} = \Delta_{32}^{-1} \Delta_{21} \Delta_{31}
\Delta_{32} .
$$
Since
$$
|V| = 2 | \Delta_1| + |\Delta_{21}| + |\Delta_{31}|
$$
and
$$
|| \Delta^{-1} V \Delta || = 2 |\Delta_{32}| + |\Delta_{21}| +
|\Delta_{31}|,
$$
we get $|\Delta_1| = |\Delta_{32}|$. Taking $V_0 \equiv 1$ we get
that $\Delta$ and $V$ are of the form as in the statement.

\smallskip

The statement (2) follows by similar considerations.
\end{proof}

\section{Constructing of $\varphi$-twisted conjugated normal form}

Let $V$ and $V'$ be words in the alphabet $\mathbb X$. If there
exists $X \in F_n$ such that $V' = \varphi(X^{-1}) V X$, we say
that $V'$ and $V$ are $\varphi$-{\it twisted conjugated}, and that
$V'$ is obtained from $V$ by $\varphi$-{\it twisted conjugation
by} $X$. Obviously, the property to be $\varphi$-twisted
conjugated is equivalence relation in the group $F_n$. Using this
definition, the Makanin's question can be reformulated as the
following: is there exists an algorithm that admits for a given
pair of elements $U$ and $V$ of a free group $F_n$ to decide if
they are $\varphi$-twisted conjugated.

In virtue of Lemma~\ref{lemma2} the equality $\varphi(\Delta^k) =
\Delta^k$ holds for any integer $k$. Therefore, if elements $U$
and $V$ are conjugated by some power of $\Delta $, they are
$\varphi$-twisted conjugated.

Below we will construct an algorithm to choose a unique
representative for each class of $\varphi$-twisted conjugated
elements. We will call this representative the {\it
$\varphi$-twisted conjugated normal form}. It will be shown that
two elements of $F_n$ are $\varphi$-twisted conjugated if and only
if their $\varphi$-twisted conjugated normal forms coincide.

If length of a reduced word $U$ in the alphabet $\mathbb X$ is
bigger than $1$ then it can be represented as a product $U \equiv
U' U''$ of two nonempty reduced words. The word $U'$ will be
referred to as an {\it initial part} of $U$ and $U''$ will be
referred as a {\it final part} of $U$. Denote by $I(U)$ the set of
all initial parts of $U$ and by $F(U)$ the set of all final parts
of $U$.

For a word $U$ and an integer $\ell$ we denote by $\varphi^{\ell}
(U)$ the word obtained from $U$ by replacing (graphically) each
letter $u$ of $U$  by its image $\varphi^{\ell}(u)$, and denote by
$U_{[\ell]}$ the word obtained after reducing of
$\varphi^{\ell}(U)$.

Let $V \equiv V' V''$ be a reduced word in the alphabet $\mathbb
X$, where $V' \in I(V)$, $V'' \in F(V)$. The word $V''_{[1]} V'$
(that represents the element $\varphi(V'') V (V'')^{-1}$ is said
to be a {\it cyclic $\varphi$-shift of a final part of $V \equiv
V' V''$} and the word $V'' V'_{[-1]}$, that represents the element
$$\varphi([\varphi^{-1} (V')]^{-1}) V \varphi^{-1}(V') = (V')^{-1}
V \varphi^{-1} (V'),$$ is said to be a {\it cyclic $\varphi$-shift
of an initial part of $V \equiv V' V''$}. If $|V'| = 1$, i.e. $V'
= x_i^{\varepsilon}$, $\varepsilon = \pm 1$, then the
corresponding $\varphi$-shift of the initial part will be referred
to as a {\it cyclic $\varphi $-shift of the initial letter}. If
$|V''| = 1$ then the corresponding $\varphi$-shift of the final
part will be referred to as a {\it cyclic $\varphi $-shift of the
final letter}.

A reduced word $V$ will be referred to as a {\it cyclically
$\varphi $-reduced} if neither a cyclic $\varphi $-shift of its
final letter nor a cyclic $\varphi $-shift of its initial letter
do not decrease length of $V$. Obviously, applying cyclic $\varphi
$-shifts of the final letter (as well as of the initial letter) to
a given word $V$, after a finite number of steps we will obtain a
cyclically $\varphi$-reduced word corresponding to $V$.

Now let us construct conjugated normal form for a word $V$.
Without loss of generality (applying, if necessary, finite number
of steps of above described $\Delta$-reducings and cyclic
$\varphi$-shifts), we can assume that $V$ is $\Delta$-reduced and
cyclically $\varphi$-reduced.

For given $V$ let ${\mathcal V}_{\Delta}$ be the set, consisting
of $V$ and all words which are conjugated to $V$ by powers of
$\Delta$ and are $\Delta$-reduced. In virtue of Lemma~\ref{lemma3}
and Lemma~\ref{lemma4}, the set ${\mathcal V}_{\Delta}$ is finite.
Applying to all elements of ${\mathcal V}_{\Delta}$ cyclic
$\varphi$-shifts of all initial parts and all final parts, we will
construct the set $({\mathcal V}_{\Delta})_{\varphi}$. Applying
$\Delta$-reducing to all element of $({\mathcal
V}_{\Delta})_{\varphi}$, we will construct the set of
$\Delta$-reduced words $(({\mathcal
V}_{\Delta})_{\varphi})_{\Delta}$. And finally, applying, if
necessary, cyclic $\varphi $-shifts of the final letter and of the
initial letter to words from $(({\mathcal V}_{\Delta})_{\varphi})
_{\Delta}$, will construct the set $D_0(V)
 = ((({\mathcal V}_{\Delta})_{\varphi})_{\Delta})_{\varphi}$ of
$\Delta$-reduced and cyclically $\varphi $-reduced words.

Recall that in the free group the set of words obtained from a
given word $V$ by cyclic shifts is finite. But the set of words
obtained from $V$ by cyclic $\varphi$-shifts can be infinite.
Indeed, it is clear from relations
$$
\left[ \varphi^{k}(V) \varphi^{k-1}(V) \cdots \varphi(V) \right] V
\left[ V^{-1} \, \varphi(V^{-1}) \cdots \varphi^{k-1}(V^{-1})
\right] = \varphi^{k}(V),
$$
for $k>0$ and
$$
\left[ \varphi^{k}(V^{-1}) \varphi^{k+1}(V^{-1}) \cdots  V^{-1}
\right] V \left[ \varphi^{-1}(V) \cdots \varphi^{k+1}(V)
\varphi^{k}(V) \right] = \varphi^{k}(V),
$$
for $k<0$, that $V$ is $\varphi$-twisted conjugated to
$\varphi^k(V)$ for any integer $k$.

For each integer $k$ we define a set $D_k(V) = D_0 (\varphi^k(V))$
and define $D(V) = \bigcup_{k \in \mathbb Z} D_k(V)$. Let us
verify that $D(V)$ is finite.

\begin{lemma} \label{lemma5}
The following equality holds:
$$
D(V) \, = \, \bigcup_{k \in \{ 0, 1, \ldots, m-1\}} D_k (V) .
$$
\end{lemma}

\begin{proof}
Let $r$, satisfying $0\leq r < m$, be such that $ k = m q + r$,
$q\in \mathbb Z$. By Lemma~\ref{lemma2},
$$
\varphi^{k} (V) \, = \, \Delta^{-q} \, \varphi^r (V) \, \Delta^q .
$$
We show that $D_k(V) = D_r (V)$, that will imply the statement.
Indeed, by the definition,  $D_k (V) = D_0 (\Delta^{-q} \varphi^r
(V) \Delta^q)$. Denote $U \equiv V_{[r]} = \varphi^r (V)$ and
consider elements from $D_0(U)  = D_0(\varphi^r(V))$. By the
definition, $D_0 (U)$ consists of words $\varphi(U'') U'$ and $U''
\varphi^{-1} (U')$ (for all pairs of initial and final parts of $U
\equiv U' U''$) to which $\Delta$-reducing and cyclic $\varphi
$-reducing are applied. To construct $D_k(V)$ we need to pass from
the word $\Delta^{-q} \varphi^r (V) \Delta^q$ to a
$\Delta$-reduced word, which, in a general case, can be different
from $U \equiv \varphi^r (V)$. But, according to the definition of
$D_0(V)$, the set of $\Delta$-reduced words constructed from
$\Delta^{-q} \varphi^r (V) \Delta^q$ coincides with the set of
$\Delta$-reduced words constructed from $\varphi^r (V)$. So,
corresponding sets of all cyclic $\varphi$-shifts of initial parts
and of final parts also coincide. Therefore, $D_k (V) = D_0
(\varphi^r(V)) = D_r (V)$.
\end{proof}

\begin{lemma}\label{lemma6}
The set $D(V)$ has the following properties: \\
(1) If $U \in D(V)$ then $D(U) = D(V)$; \\
(2) If $V$ and $W$ are $\varphi$-twisted conjugated, and each of
them is $\Delta$-reduced and cyclically $\varphi$-reduced, then
$D(V) = D(W)$.
\end{lemma}

\begin{proof}
Since $D(V)$ consists of words which are $\varphi$-twisted
conjugated in $F_n$, item (2), obviously, implies (1). Let us
prove (2). Firstly we remark that if $W$ is obtained from $V$ by a
$\varphi$-shift of an initial part or a $\varphi$-shift of a final
part, then, obviously, $D(W) = D(V)$. Now assume that $W$ and $V$
are related by $ W \, = \, U^{-1}_{[1]} V U$ for some reduced $U
\in F_n$. Since $W$ is cyclically $\varphi$-reduced, the product
$U^{-1}_{[1]} V U$ contains cancellations. Moreover, these
cancellations are either in the product $U^{-1}_{[1]} V$ or in the
product $V U$, but not in the both, since by the assumption $V$ is
$\varphi$-reduced. There are two possibilities: $V$ will be
cancelled wholly or not.

\smallskip

\underline{Case 1:} Suppose $V$ is not cancelling wholly.

\smallskip

\underline{Case 1(i):} Suppose that there are cancellations in the
product $U^{-1}_{[1]} V$.  Then $U^{-1}_{[1]}$ must be cancelled
wholly, since $W$ is assumed to be $\varphi$--reduced, and so $V
\equiv U_{[1]} V_1$. Then
$$
W \, = \, U^{-1}_{[1]} V U \, \equiv \, U^{-1}_{[1]} (U_{[1]} V_1)
U \, =  \, V_1 U
$$
and $W \equiv  V_1 U$ is obtained from $V \equiv U_{[1]} V_1$ by
the $\varphi$-shift of the initial part $U_{[1]}$. Therefore,
$D(V) = D(W)$.

\smallskip

\underline{Case 1(ii):} Suppose that there are cancellations in
the product $VU$.  Then $U$ must be cancelled wholly since $W$ is
assumed to be $\varphi$--reduced, and so $V \equiv V_1 U^{-1}$.
Then
$$
W \, = \, U^{-1}_{[1]} V U \, \equiv \, U^{-1}_{[1]} (V_1 U^{-1})
U \, = \, U^{-1}_{[1]} V_1
$$
and $W \equiv U^{-1}_{[1]} V_1$ is obtained from $V \equiv
 V_1 U^{-1}$ by the $\varphi$-shift of the final part $U^{-1}$.
Therefore, $D(V) = D(W)$.

\smallskip

\underline{Case 2:} Suppose that  $V$ is cancelling wholly. We
will use induction by length of the word $U$.

\smallskip

\underline{Case 2(i):} Suppose that $V$ is cancelling wholly in
the product $VU$ and after that there are cancellations of letters
of the remaining part of $U$ with letters of $U^{-1}_{[1]}$. Then
we can represent $U \equiv V^{-1} U_1$, therefore $U^{-1}_{[1]} =
U^{-1}_{1 \, [1]} V_{[1]}$ and
$$ U^{-1}_{[1]} V U
\, = \,  U^{-1}_{1 \, [1]} V_{[1]} V V^{-1} U_1 \, = \, U^{-1}_{1
\, [1]} V_{[1]} U_1 . $$ If $U_1$ is cancelling wholly with
$V_{[1]}$, then $V_{[1]} \equiv V_2 U_1^{-1}$ and
$$
U^{-1}_{1 \, [1]} V_{[1]} U_1 \, \equiv \, U^{-1}_{1 \, [1]} V_2
U^{-1}_1 U_1 \, = \, U^{-1}_{1 \, [1]} V_2 \, =  \, W .
$$
Remark that $W$ arises in the process of the construction of $D_1
(V)$. Indeed,
$$
D_1(V) = D_0 (V_{[1]}) = D_0 ( V_2 U_1^{-1})
$$
and $W$ is obtained from $V_2 U_1^{-1}$ by the cyclic
$\varphi$-shift of the final part. If, again, $V_{[1]}$ is
cancelling wholly with $U_1$, then the statement follows from the
induction assumption.

\smallskip

\underline{Case 2(ii):} Let there be a cancellation in the product
$U^{-1}_{[1]} V$ with $V$ cancelling wholly and remaining part of
$U^{-1}_{[1]}$ cancelling with $U$, i.e. $U^{-1}_{[1]} \equiv
U_1^{-1} V^{-1}$, $U = \varphi^{-1} (V) \varphi^{-1} (U_1)$ and
\begin{equation}
U^{-1}_{[1]} V U = U^{-1}_1 V^{-1} V V_{[-1]} U_{1 \, [-1]} \, =
\, U^{-1}_1 V_{[-1]} U_{1 \, [-1]} . \label{eq1}
\end{equation}
If $U_1^{-1}$ is cancelling wholly with $V_{[-1]}$, then $V_{[-1]}
= U_1 V_2$ and we get $U_1^{-1} V_{[-1]} U_{1 \, [-1]} = U_1^{-1}
(U_1 V_2) U_{1 \, [-1]} = V_2 U_{1 \, [-1]}$, i.e. $V = U_{1 \,
[1]} V_{2 \, [1]}$ and $W \equiv V_2 U_{1 [-1]}$. Comparing these
words we see that $W$ belongs to $D_{-1}(V)$. Indeed, by the
definition,
$$
D_{-1} (V) \, = \, D_{0} (\varphi^{-1} (V)) \, = \, D_0 ( U_1 V_2)
.
$$
Applying to $U_1 V_2$ the cyclic $\varphi$-shift of the initial
part $U_1$, we will get $W = V_2 U_{1 \, [-1]}$. The case when
$V_{[-1]}$ is cancelling wholly with $U^{-1}_1$ in (\ref{eq1})
follows from the induction assumption.
\end{proof}

\smallskip

Now we are able to complete the proof of Theorem~\ref{theorem1}.

\begin{proof}
For given $V$ we have constructed the finite set $D(V)$. From this
set we choose words of minimal length, and after that,  from such
words choose the word which is minimal in respect to the above
defined ordering on ${\mathbb F}_n$. Denote this word by
$\overline{V}$ and call it the {\it normal $\varphi$-twisted
conjugated form for $V$}.

Let $U$ and $V$ be reduced words in the free group $F_n$. Let us
construct words $\overline{U}$ and $\overline{V}$, which are
normal $\varphi$-conjugated forms for $U$ and $V$, respectively.
In virtue of Lemma~\ref{lemma5}, $\overline{U} = \overline{V}$ if
and only if elements $U$ and $V$ are $\varphi$-twisted conjugated
in $F_n$, that means, by the definition, that equation $\varphi(X)
U = V X$ is solvable in $F_n$. The proof is completed.
\end{proof}

%



\end{document}